\documentclass[12pt]{amsart}

\usepackage{amssymb}
\usepackage{color}

\usepackage[hmargin=2.5cm,bmargin=2.5cm,tmargin=3cm]{geometry}
\usepackage{graphicx}

\newtheorem{theorem}{Theorem}[section]

\theoremstyle{remark}

\newcommand{\R}{\mathbb{R}}
\providecommand{\C}{\mathbb{C}}
\renewcommand{\C}{\mathbb{C}}

\newcommand\GroupS{\mathfrak{S}}
\newcommand\Op{L}
\newcommand\Space{\mathcal{H}}
\newcommand\letter{letter} 
\newcommand\Graph{\mathsf{G}}
\newcommand\Cube{\mathsf{G}_c}
\newcommand\me{\mathsf{e}}
\newcommand\mv{\mathsf{v}}

\usepackage{bm}

\title{Eigenspaces of symmetric graphs are not typically irreducible}

\author{Gregory Berkolaiko}
\address{Department of Mathematics, Texas A\&M University, College
 Station, TX 77843-3368, USA}
\author{Wen Liu}
\address{Department of Mathematics, Lamar University, Beaumont, TX 77710, USA}

\begin{document}

\begin{abstract}
  We construct rich families of Schr\"odinger operators on symmetric
  graphs, both quantum and combinatorial, whose spectral degeneracies
  are persistently larger than the maximal dimension of an
  irreducible representation of the symmetry group.
\end{abstract}


\maketitle

\section{Introduction}

In many circumstances, it has been established that a ``typical''
Schr\"odinger-type operator $\Op$ on the Hilbert space $\Space$ has simple
spectrum.  Classical results of Uhlenbeck \cite{Uhl_ajm76} concern the Laplace
operator on compact manifolds: the spectrum can be made simple by a small
perturbation of the manifold's metric.  For quantum graphs, similar results have
been established by Friedlander \cite{Fri_ijm05} and extended by Berkolaiko and
Liu \cite{BerLiu_jmaa17}.  The role of the metric here is played by the lengths
of the graph's edges.  In the case of a combinatorial graph Laplacian, the
degeneracies in the spectrum can be lifted, for example, by changing the edge
weights or by addition of the small potential (i.e.\ a diagonal matrix).

If $\Op$ commutes with a group of unitary operators $\GroupS$, it is easy to see
that the eigenspace $E_\lambda(\Op)$ corresponding to the eigenvalue $\lambda$
is a representation of the group $\GroupS$ \cite{Wigner_grouptheory}.  It is
therefore expected\footnote{It may happen that the space $\Space$ is not rich
  enough to support some of the representations; we will see an example
  in this \letter{}.} that some of the eigenvalues will be degenerate, with the
size of the degeneracies dictated by the degrees of the group's irreducible
representations.  In analogy to the above results, it is natural to assume here
that the eigenvalues will not be more degenerate than necessary: a perturbation
with the same symmetry $\GroupS$ can ensure that every $E_\lambda(\Op)$ is
irreducible.  We will call this assumption ``generic irreducibility''.

The only rigorous result establishing generic irreducibility is due to Zelditch
\cite{Zel_aif90}, who considered the Laplace operator on finite $C^\infty$
Riemannian covers and established the positive result under the assumption that
the dimension of the manifold is greater than or equal to the maximal dimension of
irreducible representations.  To quote \cite{Zel_aif90}, ``[it] leaves
open many interesting cases of the generic irreducibility question [\ldots], in
particular, it does not touch the case of graphs''.  The purpose of this \letter{}
is to construct a rich family of graph examples (both quantum and combinatorial)
on which generic irreducibility fails.

We stress that our examples are families of Schr\"odinger operators where we
allow for perturbations not only of the metric (edge lengths or weights) but
also of the potential, as long as the prescribed symmetry is preserved.  This
contrasts the positive result of \cite{Zel_aif90} where only perturbations of
the metric were enough to resolve the degeneracy.

We will start with the quantum graphs, on which our example was originally
constructed and on which it has particularly rich structure.  We will then point
out how to translate our example to the case of combinatorial graph.  We remark
that a one-parameter family of Hubbard Hamiltonians (a model similar to a
combinatorial graph) with persistently reducible eigenstates has been previously
considered by Heilmann and Lieb \cite{HeiLie_tnyas71} (see also
\cite{Bra_pla98}).

\section{Definitions}

Let $\Graph$ be a graph with each edge $\me$ being identified with an interval
$[0,\ell_{\me}]$ of the real line. This gives us a local variable $x_{\me}$ on
the edge $e$ which can be interpreted geometrically as the distance from the
initial vertex. Which of the two end-vertices is to be considered initial is
chosen arbitrarily; the analysis is independent of this choice.

We are interested in the eigenproblem of the Schr\"odinger operator $\Op :=
-\Delta + Q$, namely
\begin{equation}
  -\frac{\partial^2}{\partial x^2} u_{\me}(x)  + Q_\me(x) u_{\me}(x) 
  = \lambda u_{\me}(x),
  \label{eq:metric_laplacian}
\end{equation}
where the potential $Q_\me(x)$ is sufficiently regular to keep the problem
self-adjoint, for instance a piecewise continuous function.  The functions $u$
are assumed to belong to the Sobolev space $H^{2}(\me)$ on each edge $\me$.  We
will impose the so-called Neumann-Kirchhoff (NK) conditions at the vertices of
the graph: we require that $u$ is continuous on the vertices, i.e.\ 
$u_{\me_{1}}(\mv)=u_{\me_{2}}(\mv)$ for each vertex $\mv$ and any two edges,
$\me_{1},\me_{2}$ incident to $\mv$, and that the current is conserved,
\begin{equation}
  \label{eq:current_cons}
  \sum_{\me\sim \mv} \frac{\partial}{\partial x} u_{\me}(\mv) = 0
  \qquad \mbox{ for all vertices } \mv,
\end{equation}
where the summation is over all edges incident to the vertex $\mv$ and the
derivative is covariant into the edge (i.e.\ if $\mv$ is the final vertex for
the edge $\me$, the corresponding term gets a minus sign). Further information can be
found in the review \cite{GnuSmi_ap06}, the textbook \cite{BerKuc_graphs} or a
recent elementary introduction \cite{Ber_prep16}, among other sources.

The symmetries we consider are induced by the graph's global
isometries, with the standard definition of the graph metric (the
length of the shortest path).  Namely, given an isometry
$s: \Graph \to \Graph$, the corresponding operator $S$ on the Hilbert
space of $L_2$ function on $\Graph$ acts as $(Su)(x) = u(s^{-1}(x))$.
It is easy to see that an isometry maps vertices to vertices,
preserving the degree and edges to edges, preserving the length.
Therefore the group of all isometries of a metric graph coincides with
the group of symmetries of the corresponding edge-weighted discrete
graph (the edge lengths become weights; symmetry transformations must
preserve weights).

As has been observed before (see, e.g. \cite{HeiLie_tnyas71}), one can
easily enrich the group of Hilbert space symmetries {\it a
  posteriori\/}, by considering all unitary operators leaving the
eigenspaces invariant.  The isometries, however, are the natural
choice of {\it a priori\/} symmetries which one would expect to
explain all degeneracies in the spectrum of a simple operator such as
Schr\"odinger.

\section{Quantum graph example}

The starting point of our considerations was an observation that the
regular tetrahedron graph (the complete graph on four vertices $K_4$
with all edge lengths equal to $a$) and no potential,
$Q_\me(x)\equiv0$, has eigenvalue $\lambda_a = (2\pi/a)^2$ with
multiplicity 4.  This is more than the maximal degree of an
irreducible representation (irrep): the symmetry group of a
tetrahedron is the symmetric group $S_4$ acting on the graph by
permuting its vertices which has irreps of dimensions $1, 1, 2, 3$ and
$3$.  Therefore, the eigenspace of $\lambda_a$ cannot be
irreducible.\footnote{It can be shown to be the sum of the standard
  and the identity representations of $S_4$.}  However, this example
is less than satisfactory, due to the paucity of the space of possible
perturbations.  The only free parameter is the length $a$ and the
entire spectrum changes trivially when all lengths are scaled by the
same factor.

It is instructive to look at the eigenspace of $\lambda_a$.  Its basis
can be chosen as follows: one eigenfunction is equal to
$\cos(2\pi x/a)$ on every edge of the graph (and is $1$ at every
vertex); three more eigenfunctions are equal to $\sin(2\pi x/a)$ on
the edges bounding one of the faces of the tetrahedron and are
identically zero on the other edges (on every vertex they are equal to
0).  The tetrahedron has four faces and one can construct four
corresponding functions, but one of them can be obtained as the sum of
the other three.  Eigenfunctions of this type give rise to many
interesting phenomena in quantum graphs, for instance to emergence of
``topological resonances'' \cite{GnuSchSmi_prl13,ColTru_prep16}, to
special terms in the zeta function of equilateral quantum graphs
\cite{HarWey_lmp17}, and to masking of the poles of the
Titchmarsh--Weyl function \cite{KuhRoh_prep16}.  And it is these
eigenfunctions which will lead us to a better example.

\begin{figure}
  \centering
  \includegraphics{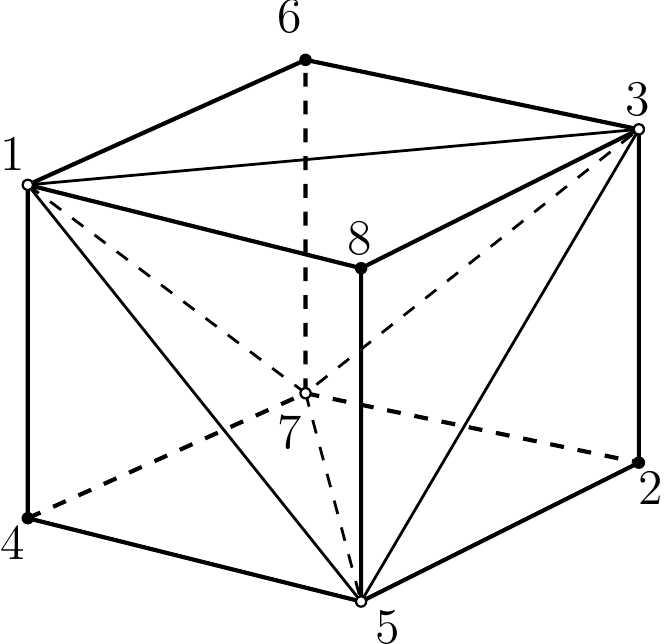}
  \caption{A graph with 8 vertices and 18 edges obtained by inscribing a tetrahedron into a cube.}
  \label{fig:cube_with_tetra}
\end{figure}

Let us inscribe a tetrahedron into a cube, see
Fig.~\ref{fig:cube_with_tetra}.  We denote this graph by $\Graph$ and
stress that the embedding of $\Graph$ into $\R^3$ is done for
visualisation reasons only; we do not assume any relation between the
length $a$ of the tetrahedron's edge and the length $b$ of the cube
edge\footnote{i.e.\ we reserve the right to change the metric on the
  edges; equivalently, we can allow the edges to be curved.}.  The
resulting graph has the symmetry of the tetrahedron.  In fact, we will
allow any potential $Q_a(x)$ and $Q_b(x)$ on the edges of length $a$
and $b$, as long as the symmetry of the graph is preserved.  This
condition only restricts $Q_a(x)$ to be even and forces the
orientation of the cube edges to be chosen consistently: all edges
oriented from the odd-numbered vertex (where the tetrahedron edges are
incident) to the even-numbered vertex.  The orientation of an edge
serves to prescribe how the potential is placed on the edge.

This graph turns out to have eigenvalues of multiplicity at
least 5, which we will demonstrate by constructing the
eigenfunctions.  We summarize this discussion as a theorem.

\begin{theorem}
  Consider the graph $\Gamma$ depicted in
  Fig.~\ref{fig:cube_with_tetra}.  Let all tetrahedron edges
  (i.e. those which connect odd-numbered vertices) have length $a$ and
  support potential $Q_a(x)$ which is even, $Q_a(a-x) = Q_a(x)$.  Let
  all cube edges (i.e. those that connect an odd-numbered vertex to an
  even-numbered one) have length $b$ and support potential $Q_b(x)$ of
  arbitrary form which is placed so that $x=0$ corresponds to the
  odd-numbered vertex and $x=b$ corresponds to the even-numbered
  vertex.

  Then the symmetry group of $\Gamma$ is the symmetric group $S_4$
  which has irreducible representations of degrees $1, 1, 2, 3$ and
  $3$, yet for any choice of $a$, $b$, $Q_a$ and $Q_b$, the graph
  $\Gamma$ has infinitely many eigenvalues of multiplicity at least
  5.  The corresponding eigenspaces must therefore be reducible.
\end{theorem}

\begin{figure}
  \centering
  \includegraphics{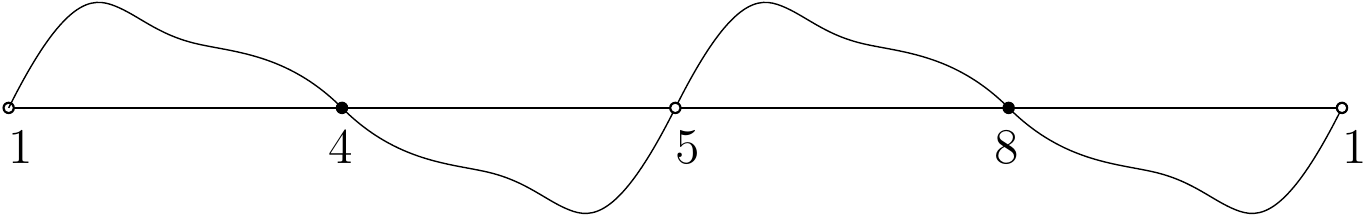}
  \caption{An eigenfunction of the Dirichlet problem \eqref{eq:Dirichlet_edge}
    replicated along the cycle $1$, $4$, $5$, $8$, $1$ to create a continuously
    differentiable function along the cycle.}
  \label{fig:eigenfunction}
\end{figure}

\begin{proof}
  Let $\lambda$ be an eigenvalue of the Dirichlet problem on the
  interval of length $b$, namely
  \begin{equation}
    \label{eq:Dirichlet_edge}
    -\psi''(x) + Q_b(x) \psi(x) = \lambda \psi(x), \qquad \psi(0) = \psi(b) = 0,
  \end{equation}
  with $\psi(x)$ the corresponding eigenfunction.  Take a cycle on the
  graph $\Graph$ consisting of the cube edges only, for example, the
  cycle on the vertices $1$, $4$, $5$, $8$.  Place the function $\psi$
  on the edges of the cycle so that the derivative is continuous along
  the cycle.  Namely, we place $\psi(x)$ on the edge $(1,4)$,
  $-\psi(b-x)$ on the edge $(4,5)$, then again $\psi(x)$ on $(5,8)$
  and so on, see Fig.~\ref{fig:eigenfunction}.  We extend the function
  by 0 to the rest of the graph.  It is easy to see that the resulting
  function is continuous on the entire graph and satisfies condition
  \eqref{eq:current_cons} on the vertices.  By virtue of equation
  \eqref{eq:Dirichlet_edge} it also satisfies the Schr\"odinger
  equation on the entire graph with the eigenvalue $\lambda$.

  Obviously, we can repeat this process for \emph{every} cycle
  consisting of edges of length $b$.  The number of the linearly
  independent functions that can be produced in this way is equal to
  the number of the linearly independent cycles on the cube subgraph
  of $\Graph$.  This, in turn, is given by the first Betti number of
  the subgraph, $\beta = E - V + 1$, which for the cube is equal to
  $5$.  Informally, the boundaries of $5$ faces of the cube are
  independent, while the sixth one is given by their sum.

  To summarize, each eigenvalue $\lambda$ of \eqref{eq:Dirichlet_edge}
  is an eigenvalue of multiplicity at least 5 of the graph $\Graph$
  and its eigenspace is a reducible representation of the group of
  symmetries of $\Graph$.
\end{proof}

\section{Variations on the example and its analysis}

We note that the role played by the tetrahedron subgraph of the graph $\Graph$ is a
very limited one: it serves to restrict the symmetry group of the resulting
graph and adds the freedom of choosing the metric on its edges (equivalently,
their length) and the potential.  We can dispense with this subgraph altogether
and consider the cube graph $\Cube$ with the odd-numbered vertices distinguished
from even-numbered.  This can be done by choosing the potential $Q_b(x)$ which
is not even, i.e.\ $Q_b(b-x) \not\equiv Q_b(x)$, or by changing the vertex
conditions at even-numbered vertices to $\delta$-type with non-zero parameter.
We will assume the former method is used and not dwell on the latter.

It is also interesting to note that it is not necessary to restrict the symmetry
of the cube: the group of cube's symmetries, the \emph{full octahedral group},
has representations of degrees up to 3 and therefore the persistent eigenspace
of dimension 5 still provides a valid counter-example to the generic
irreducibility conjecture.  However, this further restricts the space of
available perturbations and makes the forthcoming analysis unwieldy due to the
large symmetry group.

We now consider the graph $\Cube$ with the symmetry group $S_4$ and identify the
decomposition of the domain of the Schr\"odinger operator into the subspaces
corresponding to the irreducible representations of $S_4$.  More precisely, 
denote by $\Space(\Cube)$ the functions on the edges of $\Cube$ that belong to
the Sobolev space $H^2$ on every edge and satisfy continuity condition and
condition \eqref{eq:current_cons} on the vertices of the graph.  For a
representation $\rho$, let $M^\rho_g$ denote the matrix corresponding to the
group element $g\in S_4$.  We
are looking for tuples $(\psi_1, \ldots, \psi_d)^T$ of functions from
$\Space(\Cube)$ which satisfy the intertwining condition
\begin{equation}
  \label{eq:intertwine}
  \begin{pmatrix}
    \psi_1(gx) \\ \cdots \\ \psi_d(gx)
  \end{pmatrix}
  = M^\rho_g 
  \begin{pmatrix}
    \psi_1(x) \\ \cdots \\ \psi_d(x)
  \end{pmatrix}.
\end{equation}
We will call the functions satisfying this condition the \emph{equivariant
  functions} for the representation $\rho$; the subspace of $\Space(\Cube)$
spanned by them is called the \emph{isotypic component} of the representation
$\rho$ and denoted by $\Space_\rho(\Cube)$.  Once the space $\Space_\rho(\Cube)$
identified, one can restrict the operator $\Op$ to this space and unitarily
reduce it to a simpler problem.  This procedure, pioneered on quantum graphs by
Band, Parzanchevski and Ben-Shach \cite{BanParBen_jpa09,ParBan_jga10} for their
study of isospectrality is called ``quotient graph construction''.  We refer the
interested reader to these papers as well as to the forthcoming work
\cite{BanBerJoyLiu_prep17} where it is formalized in terms of the scattering
matrices.

\begin{figure}
  \centering
  \includegraphics{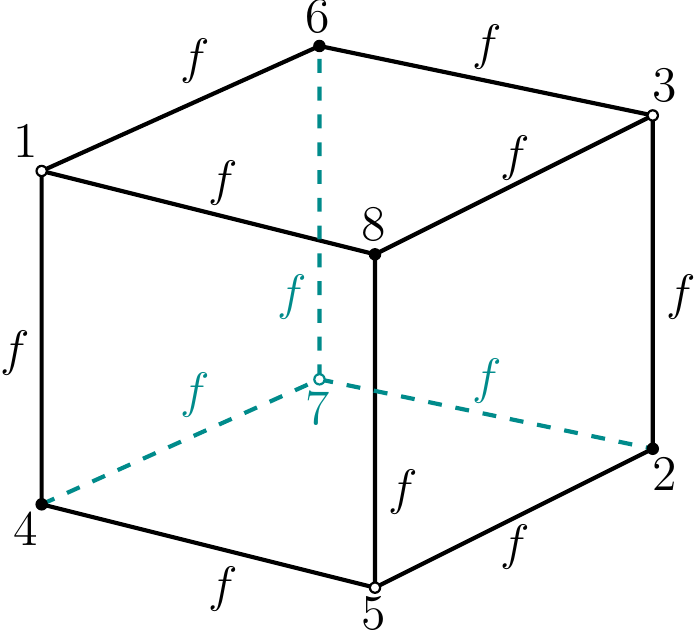}
  \caption{A representation of the functions from $\Space(\Cube)$ that transform
    according to the identity representation \eqref{eq:identity_rep}.}
  \label{fig:isotypic_i}
\end{figure}

As mentioned above,
$S_4$ has 5 irreps of degrees $1, 1, 2, 3$ and $3$.  The \emph{identity representation} maps every $g \in S_4$ to
multiplication by one,
\begin{equation}
  \label{eq:identity_rep}
  R_{i} = \Big\{ (1\,3) \mapsto (1),\quad (1\,5) \mapsto (1),\quad (1\,7) \mapsto (1) \Big\}.
\end{equation}
In this notation, for each $g$ from a set of generators of $S_4$ (here
we took $(1\,3)$, $(1\,5)$ and $(1\,7)$; we remind that the symmetry
transformations act as permutations on the odd-numbered vertices) we
specify a $1\times1$ matrix (in this case, multiplictation by $1$).
The intertwining condition becomes $\psi(gx) = \psi(x)$ for every
$g\in S_4$ which is satisfied by $\psi(x)$ which are equal to the same
function, which we denote by $f$, on every edge oriented from odd to
even-numbered vertex, see Fig.~\ref{fig:isotypic_i}.  Not every $f$ is
admissible: condition \eqref{eq:current_cons} at vertex $1$, for
example, becomes the condition $3f'(v_1) = 0$; the same at vertex $8$.
The continuity at each vertex is, of course, automatic.  The space of
functions of this form on $\Cube$ is the isotypic component of $R_{i}$
denoted by $\Space_i(\Cube)$.  For the function from $\Space_i(\Cube)$
to be an eigenfunction of our Schr\"odinger operator $\Op$ on the
graph $\Cube$, $f$ must be an eigenfunction of the Neumann problem
\begin{equation}
  \label{eq:Neumann_edge}
  -\psi''(x) + Q_b(x) \psi(x) = \lambda \psi(x), \qquad \psi'(0) = \psi'(b) = 0
\end{equation}
on the interval $[0,b]$.  The spectrum of \eqref{eq:Neumann_edge} coincides with
the spectrum of $\Op$ restricted to the space $\Space_i(\Cube)$ (actually, the
corresponding operators are unitarily equivalent).

The \emph{sign representation} maps every $g \in S_4$ into multiplication by the
sign of the permutation $g$,
\begin{equation}
  \label{eq:sign_rep}
  R_{s} = \left\{ (1\,3) \mapsto (-1),\quad (1\,5) \mapsto (-1),\quad (1\,7)
    \mapsto (-1) \right\}.
\end{equation}
It is easy to see that the space of functions satisfying \eqref{eq:intertwine}
with representation $R_s$ is the trivial space, $\Space_s(\Cube) = \{0\}$.
Indeed, taking for example the edge $(1,4)$, we observe that it is fixed by the
reflection $(5\,7)$, therefore the component of $\psi$ on this edge must satisfy
$\psi_{(1,4)}(x) = - \psi_{(1,4)}(x)$, and therefore $\psi_{(1,4)} \equiv 0$.
It is easy to verify that each edge is similarly fixed by some transposition, so
$\psi$ must be 0 on every edge.  Therefore, the operator $L$ has no
eigenvalues corresponding to the sign representation!

\begin{figure}
  \centering
  \includegraphics{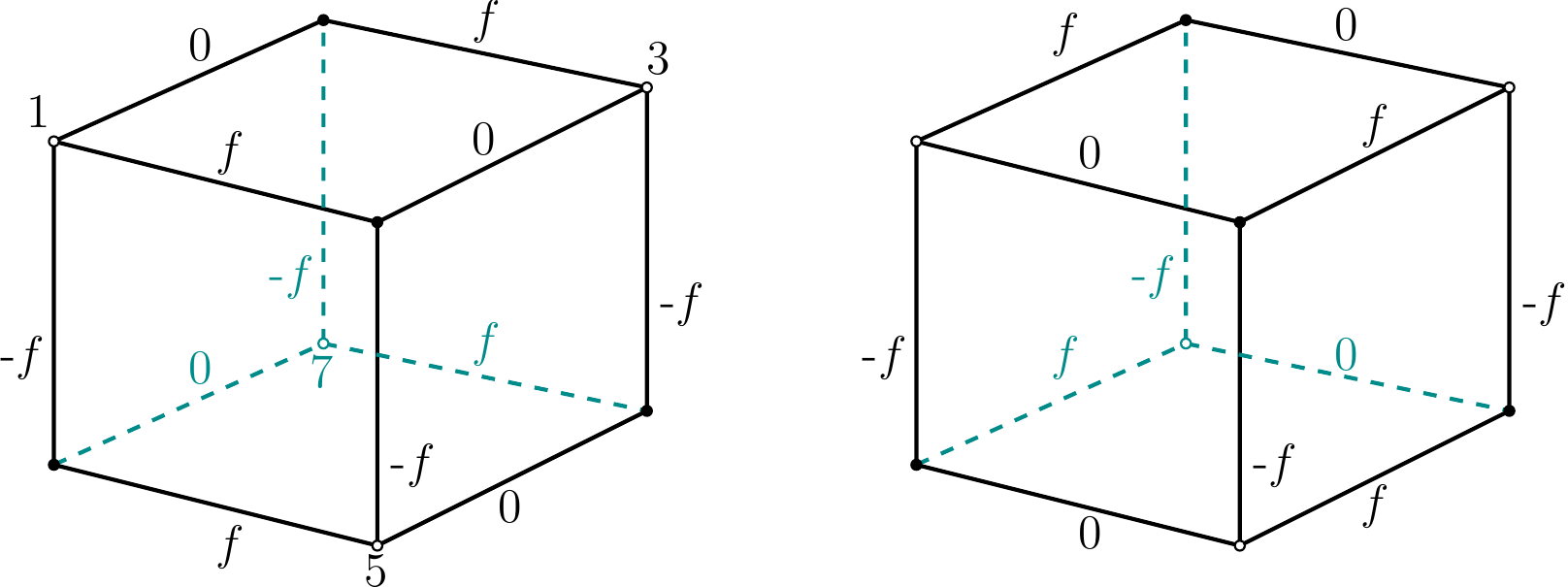}
  \caption{The structure of the functions on the graph $\Cube$ that transform
    according to the representation $R_{2d}$, equation~\eqref{eq:2d_rep}.}
  \label{fig:isotypic_2d}
\end{figure}

The next representation is the \emph{irreducible representation of degree 2}, given
in the matrix form by
\begin{equation}
  \label{eq:2d_rep}
  R_{2d} = \left\{ (1\,3) \mapsto 
    \begin{pmatrix}
      0 & 1\\
      1 & 0
    \end{pmatrix},
    \quad (1\,5) \mapsto 
    \begin{pmatrix}
      -1 & 0 \\
      -1 & 1
    \end{pmatrix},
    \quad (1\,7) \mapsto
    \begin{pmatrix}
      1 & -1 \\
      0 & -1
    \end{pmatrix}
  \right\}.
\end{equation}
Note that the matrices in this presentation are not unitary, but can be made so
using a change of basis.  However, with matrices in this form, the pairs of
functions from $\Space(\Cube)$ transforming according to this representation
have especially simple form, depicted in Fig.~\ref{fig:isotypic_2d}.  It is
immediate from the figure that for $f$ to be admissible, it must satisfy
Dirichlet problem \eqref{eq:Dirichlet_edge}.  Each admissible $f$ gives rise to
a two-dimensional eigenspace of $\Op$.

\begin{figure}
  \centering
  \includegraphics{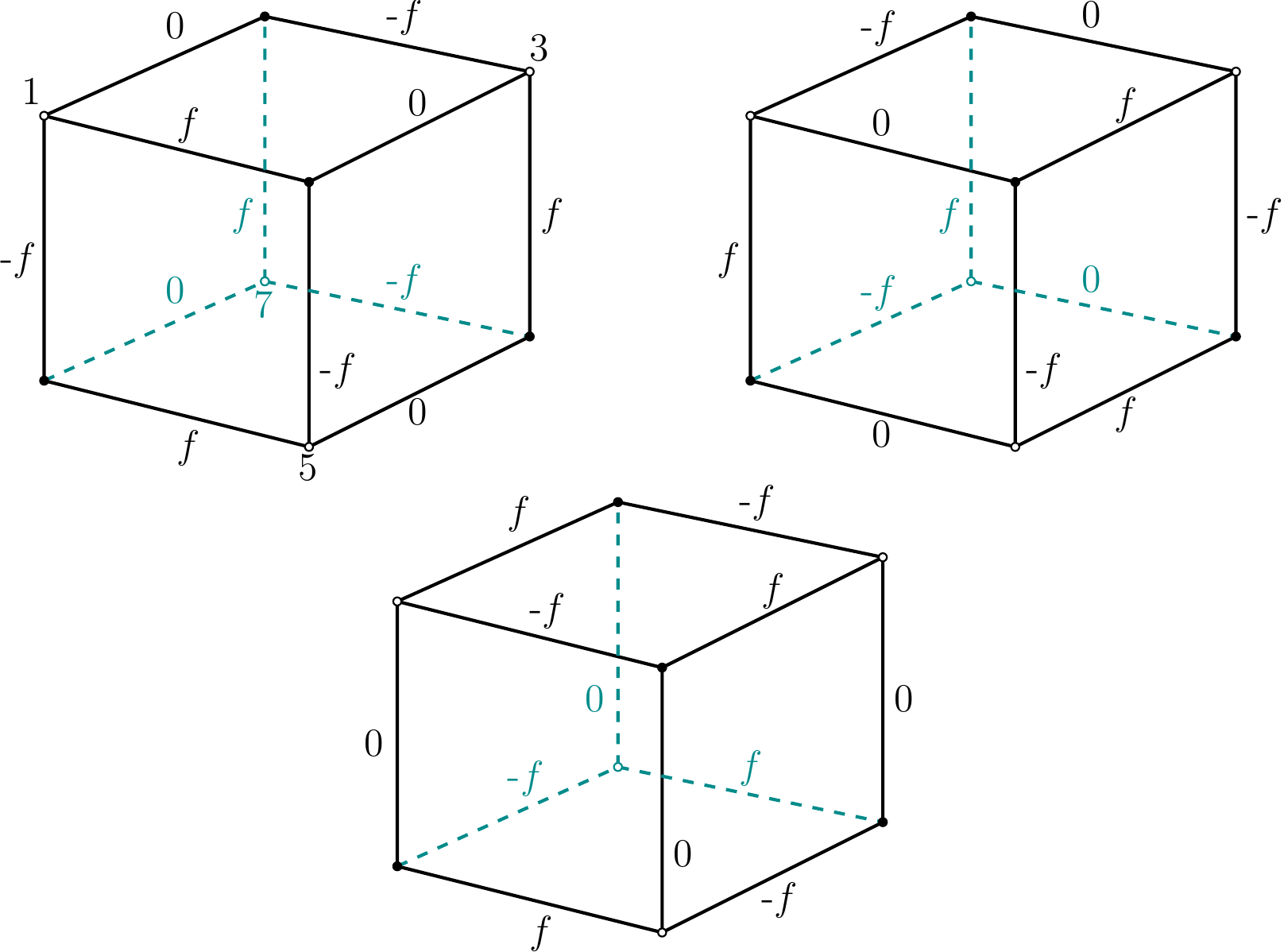}
  \caption{The structure of the functions on the graph $\Cube$ that transform
    according to the representation $R_{3d,1}$, equation~\eqref{eq:3d1_rep}.}
  \label{fig:isotypic_3d1}
\end{figure}

The \emph{standard representation} of $S_4$ is a representation of degree 3,
given by
\begin{equation}
  \label{eq:3d1_rep}
  R_{3d,1} = \left\{ (1\,3) \mapsto 
    \begin{pmatrix}
      0 & 1 & 0\\
      1 & 0 & 0\\
      0 & 0  & -1
    \end{pmatrix},
    \quad (1\,5) \mapsto 
    \begin{pmatrix}
      -1 & 0 & 0 \\
      0 & 0 & 1 \\
      0 & 1 & 0
   \end{pmatrix},
    \quad (1\,7) \mapsto
    \begin{pmatrix}
      0 & 0 & 1 \\
      0 & -1 & 0 \\
      1 & 0 & 0
    \end{pmatrix}
  \right\}.
\end{equation}
The triple of functions transforming according to this representation is
schematically represented in Fig.~\ref{fig:isotypic_3d1}. For $f$ to be
admissible, it must again satisfy Dirichlet problem \eqref{eq:Dirichlet_edge}.
We remark that the first two equivariant functions are similar in structure to
the equivariant functions we found for the representation $R_{2d}$ but differ
from them in sign distribution.

\begin{figure}
  \centering
  \includegraphics{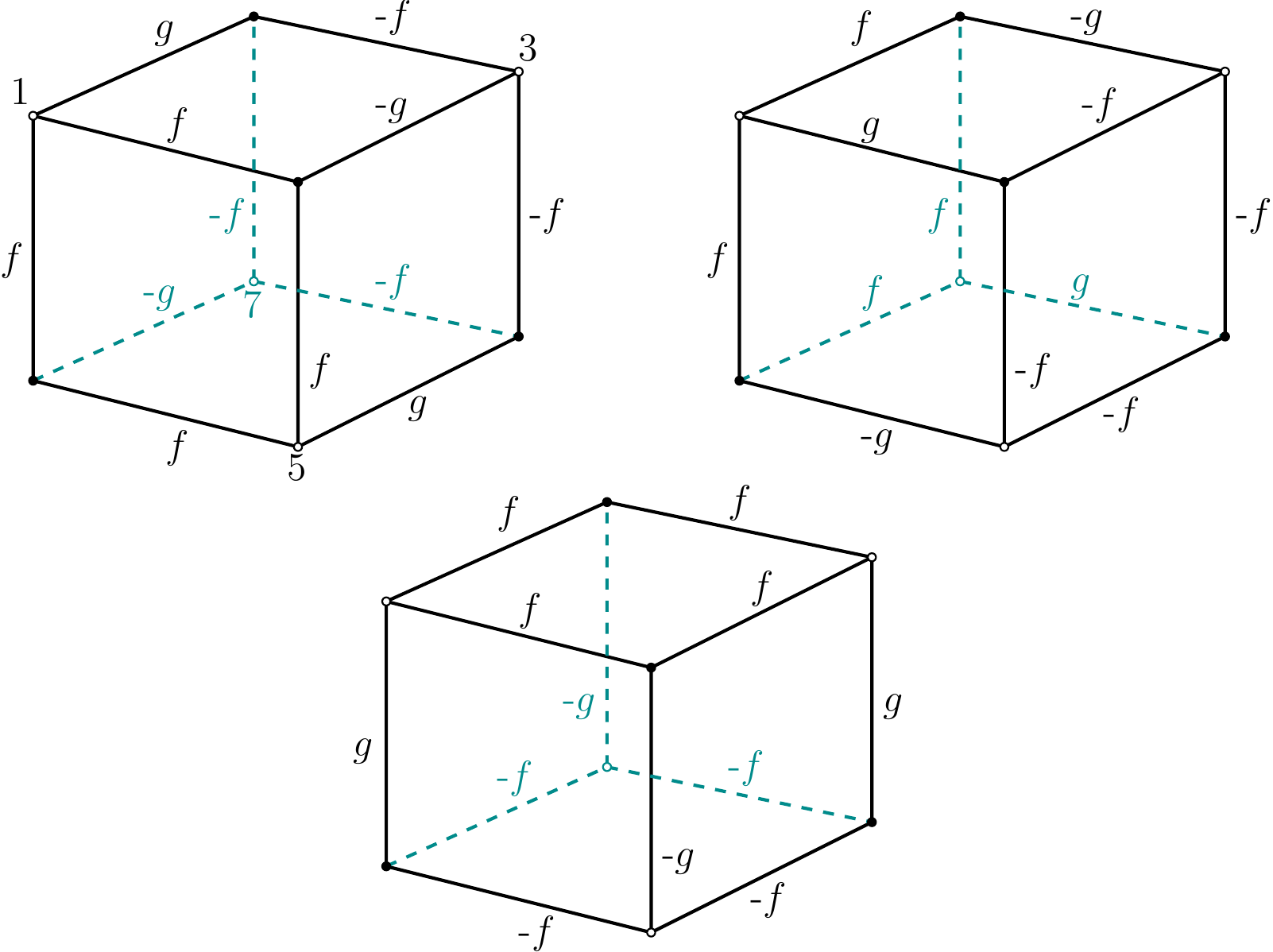}
  \caption{The structure of the functions on the graph $\Cube$ that transform
    according to the representation $R_{3d,2}$, equation~\eqref{eq:3d2_rep}.}
  \label{fig:isotypic_3d2}
\end{figure}

Finally, the last irreducible representation is the \emph{product of the
  standard and sign representations}.  It has degree 3 and is given by
\begin{equation}
  \label{eq:3d2_rep}
  R_{3d,2} = \left\{ (1\,3) \mapsto 
    \begin{pmatrix}
      0 & -1 & 0\\
      -1 & 0 & 0\\
      0 & 0  & 1
    \end{pmatrix},
    \quad (1\,5) \mapsto 
    \begin{pmatrix}
      1 & 0 & 0 \\
      0 & 0 & -1 \\
      0 & -1 & 0
   \end{pmatrix},
    \quad (1\,7) \mapsto
    \begin{pmatrix}
      0 & 0 & -1 \\
      0 & 1 & 0 \\
      -1 & 0 & 0
    \end{pmatrix}
  \right\}.
\end{equation}
The triple of functions transforming according to this representation is
schematically represented in Fig.~\ref{fig:isotypic_3d2}.  Assume $x=0$ at an
odd-numbered vertex and $x=b$ at an even-numbered vertex.  Then, for $f$ and $g$
to be admissible, they must satisfy the following problem
\begin{align}
  \label{eq:fg_on_edge}
    -f''(x) + Q_b(x) f(x) &= \lambda f(x), &
    f(0) &= g(0), &
   2f'(0) + g'(0) &= 0, \\
    -g''(x) + Q_b(x) g(x) &= \lambda g(x), &
    f(b) &= -g(b), &
    2f'(b) - g'(b) &= 0.
\end{align}
This problem is self-adjoint in an appropriately weighted $L^2\times L^2$ space.

To summarize, the eigenvalues $\lambda$ of the Dirichlet problem
\eqref{eq:Dirichlet_edge} on a single edge are also present in the spectrum of
the operator $\Op$ on the graph $\Cube$ with multiplicity at least 5.  Their
subspaces reduce into the direct sum of the degree-two and the standard
representations of the symmetry group of the underlying graph.

\section{An example of a combinatorial graph with reducible eigenspaces}

\begin{figure}
  \centering
  \includegraphics{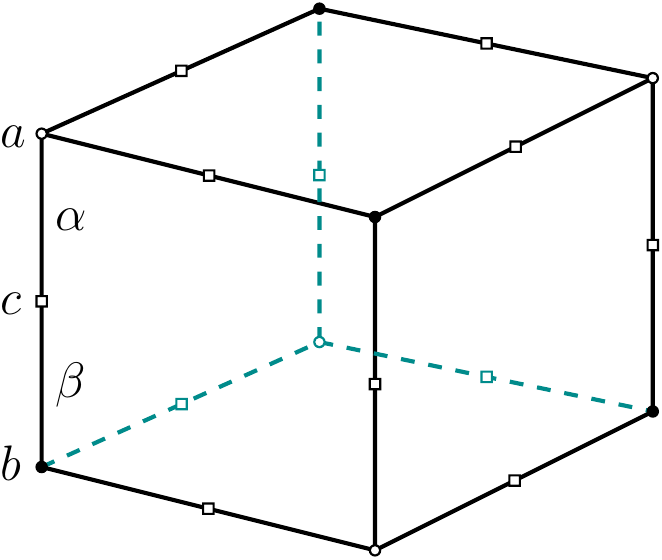}
  \caption{An example of a combinatorial graph whose Schr\"odinger operator has
    a reducible eigenspace for all values of the vertex potential $a$, $b$ and
    $c$ and for all coupling weights $\alpha$ and $\beta$.}
  \label{fig:discrete}
\end{figure}

It is easy to construct an example of a combinatorial graph with reducible
eigenspaces by analogy with the quantum graph $\Cube$.  The simplest such
example is shown in Fig.~\ref{fig:discrete}.  The corresponding Schr\"odinger
operator is a $20\times20$ self-adjoint matrix with 5 real parameters.  Namely,
we can choose the potential at 3 types of vertices ($a$, $b$ and $c$ in the
picture) and the coupling weights corresponding to two types of edges ($\alpha$
and $\beta$ in the picture).  One can construct eigenvectors by choosing a face
and placing alternating $\pm1$ on the $c$-type vertices around that face; all
other entries of the vector are zero.  It is easy to check that it is indeed an
eigenvector with the eigenvalue equal to the potential $c$.  Similarly to above,
there are 5 such independent eigenvectors.

It is straightforward to deduce the quotient eigenvalue problems as was done for
the quantum graph $\Cube$ above.  The result is analogous and we leave the
details to reader; the theory of constructing quotient combinatorial graphs will
be formalized in \cite{BanBerJoyLiu_prep17}.

\section{Concluding remarks}

Since the operators we considered have real coefficients, they are also
symmetric with respect to complex conjugation.  The choice of the field ($\R$ or
$\C$) over which the representation is irreducible can play an important role
(see the remarks in the end of \cite[Sec. 1b]{Zel_aif90}; see also
\cite{BerCom_prep14} for a different example).  However, for the symmetry group
in our example, the irreducible representations of the symmetry group over real
numbers and over complex numbers coincide.

It is unclear at the moment if it is possible to predict (without direct
computation) that the quotient graphs by $R_{2d}$ and by $R_{3d,1}$ will have
coinciding spectra.

One may speculate that the large multiplicities in the examples we
constructed may be viewed as traces of larger symmetry groups of
2-dimensional graph-like manifolds that were shrunk to the graph limit
(see \cite{Post_book12} and references therein).  As a starting point for
this process one may take the celebrated Klein's quartic, a compact
Riemann surface in the shape of a tetrahedron with the highest
possible order (namely, 168) automorphism group for its genus \cite{KleinsQuartic_book}.

However, we feel that the central role in this example is played not
be 1-dimensionality of the edges, but by vertices: quantum graphs are
not 1-dimensional manifolds as they singular at the vertices.  As a
consequence, the unique continuation principle fails on graphs.  In
fact, one can create a host of similar examples by modifying a graph
with a large symmetry group with a choice of few rank one
perturbations at the vertices (for example, changing NK conditions to
Dirichlet).  Few well-placed modifications can completely break the
symmetry yet each rank-one perturbation will split off only a single
eigenvalue from each degenerate group.

\section*{Acknowledgement}

Both authors were partially supported by the NSF grant DMS-1410657.
The authors are grateful to Ram Band for pointing out the simple
geometric meaning of our originally clumsy example.  Ram Band, Chris
Joyner and Peter Kuchment have patiently listened to our explanations
at various stages of development and helped us along with
criticism and encouragement.  We thank Jan Segert for asking difficult
questions and for pointing out the possible beautiful connection to
the Klein's quartic.



\bibliographystyle{plain}
\bibliography{bk_bibl,additional}

\end{document}